\documentclass[11pt]{amsart}
\usepackage{amsmath,amssymb,amstext}

\makeatletter \@addtoreset{equation}{section}

\makeatother \makeatletter

\newcommand{\N}{\ensuremath{\mathbb{N}}}

\newcommand{\R}{\ensuremath{\mathbb{R}}}
\newcommand{\C}{\ensuremath{\mathbb{C}}}

\newcommand{\T}{\ensuremath{(T_t)_{t\geq 0}}}

\newcommand{\supp}{\operatorname{supp}}

\newcommand{\norm}[1]{{\|#1\|}}

\newcommand{\abs}[1]{\lvert#1\rvert}

\DeclareMathOperator{\dens}{dens}
\DeclareMathOperator{\Dens}{Dens}
\newcommand{\ldens}{\underline{\dens}}
\newcommand{\lDens}{\underline{\Dens}}
\newcommand{\spa}{\operatorname{span}}

\newtheorem{theorem}{Theorem}[section]
\newtheorem{corollary}[theorem]{Corollary}

\newtheorem{lemma}[theorem]{Lemma}
\newtheorem{proposition}[theorem]{Proposition}

\theoremstyle{definition}
\newtheorem{definition}[theorem]{Definition}
\newtheorem{example}[theorem]{Example}

\newtheorem{remark}[theorem]{Remark}

\begin{document}
   \title{Frequently hypercyclic semigroups}
    \author[E. M. Mangino]{Elisabetta M. Mangino$^1$}
    \address{$^1$ Dipartimento di Matematica ``Ennio di Giorgi''\\Universit\`{a} del Salento\\  I-73100 Lecce, Italy}
    \email{elisabetta.mangino@unisalento.it}
    \author[A.Peris]{Alfredo Peris$^2$}
    \address{$^2$ IUMPA, Departament de Matem\`{a}tica Aplicada, Edifici 7A\\
     Universitat Polit\`ecnica de Val\`encia
     \\ E-46022 Val\`{e}ncia, Spain}
    \email{aperis@mat.upv.es}

\begin{abstract}
We study frequent hypercyclicity in the context of strongly
continuous semigroups of operators. More precisely, we give a
criterion (sufficient condition) for a semigroup to be frequently
hypercyclic, whose formulation depends on the Pettis integral.  This
criterion can be verified in certain cases in terms of the
infinitesimal generator of semigroup. Applications are given for
semigroups generated by Ornstein-Uhlenbeck operators, and especially
for translation semigroups on weighted spaces of $p$-integrable
functions, or continuous functions that, multiplied by the weight,
vanish at infinity.
\end{abstract}

\maketitle

\bigskip



\section{Introduction}

The hypercyclic  behaviour of strongly continuous one parameter semigroups was studied in a
systematic way for the first time in the paper by Desch, Schappacher, and Webb \cite{dsw}.
They gave a sufficient condition for
hypercyclicity of a semigroup based on the analysis of the point spectrum of the generator
of the semigroup. Moreover they  characterized hypercyclic translation semigroups defined on
weighted spaces of continuous or integrable functions on the real line.
However, during the last years it was shown that hypercyclicity appears in $C_0$-semigroups
associated to  ``birth and death'' equations for cell populations, transport equations, first order
partial differential equations and diffusion operators as  Ornstein-Uhlenbeck operators
(see \cite{cm1} for a survey on the subject. Further references on hypercyclic semigroups and
related properties are, e.g.,
\cite{banasiak2005birth,banasiak_moszynski2005a,bermudez_bonilla_conejero_peris2005hypercyclic,%
bermudez_bonilla_emamirad2005chaotic,bernal-gonzalez_grosse-erdmann2007existence,conejero_peris2005linear,costakis_peris2002hypercyclic,%
desch_schappacher2005on,el-mourchid_metafune_rhandi_voigt2008on,ji_weber2010dynamics,kalmes2006on,kalmes2007hypercyclic,%
kalmes2009hypercyclic,ma1,rudnicki2004chaos,takeo2005chaos,weber2009tensor}).

We recall that, if $X$ is a separable infinite-dimensional Banach space,
    a $C_0$-semigroup $(T_t)_{t\geq 0}$ of linear and continuous operators on $X$ is said to be \emph{hypercyclic}
    if there exists $x\in X$ (called hypercyclic vector for the semigroup) such that the set $\{T_tx\,:\,t\geq 0\}$ is dense in $X$.
    An element $x\in X$ is said to be a \emph{periodic point} for the semigroup if there exist $t>0$ such that $T_tx=x$.
    A semigroup $(T_t)_{t\geq 0}$ is called \emph{chaotic} if it is hypercyclic and the set of periodic points is dense in $X$.\par\medskip

In \cite{conejero_muller_peris2007hypercyclic}, the second author,
in collaboration with A. Conejero and V. M\"uller, proved that if
$x\in X$ is a hypercyclic vector for $\T$ then  for every $t>0$ the
set $\{T_{nt}x\,:\,n\in \N\}$ is dense in $X$, i.e. $x$ is a
hypercyclic vector for each single operator $T_t$, $t>0$. In
particular, hypercyclicity is inherited by discrete subsemigroups.
However, this is not the case in general if we change the index set
\cite{conejero_peris2009hypercyclic}, or if we consider the chaos
property \cite{bayart_bermudez2009semigroups}.

    Motivated by Birkhoff's ergodic
        Theorem, Bayart and Grivaux  introduced the notion of
        frequent hypercyclic operators \cite{bayart_grivaux2006frequently}  (see  \cite{bama} and the references quoted therein, see also
        \cite{BGE1,BGE,grosse-erdmann_peris2005frequently}),
        trying to quantify the frequency with which an orbit meets the open sets.
        This concept was extended to $C_0$-semigroups in \cite{badea_grivaux2007unimodular}.
        We recall that  the \emph{lower density} of a measurable set $M\subset
        \mathbb{R}_+$ is defined by

        \[
        \lDens (M):=\liminf_{N\rightarrow \infty} \mu(M\cap [0,N])/N,
        \]
        where $\mu$ is the
        Lebesgue measure on $\mathbb{R}_+$.
              A $C_0$-semigroup $\T$ is said to be
        \emph{frequently hypercyclic} if there exists $x\in X$ such that
         $\lDens (\{t\in \mathbb{R}_+:T_tx\in
        U\})>0$ for any non-empty open set $U\subset X$.

         If the lower density of a set $A\subset
        \mathbb{N}$ is defined by

        \[
        \ldens (A):=\liminf_{N\rightarrow \infty} \#\{n\leq N:n\in
        A\}/N,
        \]
         an operator $T\in L(X)$ is said to be \emph{frequently
        hypercyclic} if there exists $x\in X$ (called frequently hypercyclic vector) such that,
        for every non-empty open subset $U\subset X$, the set $\{n\in\mathbb{N}:T^nx\in
        U\}$ has positive lower density.
         In \cite{conejero_muller_peris2007hypercyclic} it was proved that if $x\in X$ is a frequently hypercyclic
         vector for $\T$, then for every $t>0$ the $x$ is a frequently hypercyclic vector for the operator $T_t$.

In \cite{BGE1,BGE}, Bonilla and Grosse-Erdmann  improve a result of
Bayart and Grivaux and show the following Frequent Hypercyclicity
Criterion for operators (See also \cite{griv} for a probabilistic
criterion).

\begin{proposition}\label{pr:frhyop} Let $T$ be a continuous operator on a separable Banach space $X$. If
there  exist $X_0\subseteq X$, dense subset, and a map $S:X_0\rightarrow X_0$  satisfying:
\begin{itemize}
\item[{\rm (i)}] $TSx=x$, for all $x\in X_0$;
\item[{\rm (ii)}] $\sum_{n=1}^\infty T^nx$ is unconditionally convergent  for all $x\in X_0$;
\item[{\rm (iii)}] $\sum_{n=1}^\infty S^nx$ is unconditionally convergent  for all $x\in X_0$.
\end{itemize}
then $T$ is frequently hypercyclic.
\end{proposition}

The aim  of the present paper is to give a continuous version of the
frequent hypercyclicity criterion. The unconditional convergence of
the series in Proposition \ref{pr:frhyop} will be replaced by the
Pettis integrability of orbits under the semigroup. Thanks to this
criterion we will show that, e.g., the well-known Desch-Schappacher
criterion for chaotic semigroups (see \cite{dsw}) is actually a
condition for frequent hypercyclicity. Moreover we prove that
chaotic translation semigroups on weighted spaces of  integrable functions defined on $[0,+\infty[$ are
frequently hypercyclic. We give a necessary condition on the weight
for frequently hypercyclicity. Since several properties of the
Pettis integral are used in the proofs, for the convenience of the
reader we recall in an appendix the main definitions and basic
results.

\section{Frequent Hypercyclicity Criterion for semigroups}

\begin{proposition} Let $\T$ be a $C_0$-semigroup on a separable Banach space $X$.
Then the following conditions are equivalent:
\begin{itemize}
\item[{\rm (i)}] $\T$ is frequently hypercyclic,
\item[{\rm (ii)}] for every $t>0$ the operator $T_t$ is frequently hypercyclic.
\item[{\rm (iii)}] there exists $t>0$ such that the operator $T_t$ is frequently hypercyclic.
\end{itemize}
\end{proposition}

\proof The implication (i) $\Rightarrow$ (ii) was proved in \cite{conejero_muller_peris2007hypercyclic}.
It remains to prove that (iii) implies (i).
Assume  w.l.o.g. that $t=1$ and let $x$ be a frequently hypercyclic vector for $T_1$.
Let $y\in X$, $U,V$  $0$-neighbourhoods such that $V+V\subseteq U$.
By the strong continuity of $\T$, there exists $0<\delta<1$ such that $T_sy-y\in V$ for every
$s\in [0,\delta]$. Moreover, by the local equicontinuity of $\T$, there exists  a $0$-neighbourhood $V'$
such that $T_s(V')\subseteq V$ for every $s\in [0,\delta]$.

By assumption,

\[
\ldens \{ n\in\N \ : \  T^nx\in y+V'\}>0.
\]

If $T^nx\in y+V'$, then for every $t\in [n,n+\delta]$

\[
T_tx-y=T_{t-n}(T_n x -y) + T_{t-n}y-y \in T_{t-n}(V')+ V\subseteq V+V\subseteq U.
\]

Thus, for every $N\in \N$

\[
\frac{\mu\{ t\leq N \ : \  T_tx \in U\}}{N} \geq \delta \frac{\#
\{n\in \N \ : \  n\leq N,\ T_nx\in y+V'\}}{N},
\]
hence

\[
\liminf_{N\to\infty} \frac{\mu\{ t\leq N \ : \  T_tx \in U\}}{N}
\geq \delta \liminf_{N\to\infty}\frac{\# \{n\in \N \ : \  n\leq N,\
T_nx\in y+V'\}}{N}>0. \qquad \qed
\]

\begin{theorem}\label{fhcriterion} Let $\T$ be a $C_0$-semigroup on a separable
Banach space $X$ such that there exist $X_0\subseteq X$, dense subset, and maps
$S_t:X_0\rightarrow X$, $t>0$, satisfying:
\begin{itemize}
\item[{\rm (i)}] $T_tS_tx=x$, $T_tS_rx=S_{r-t}x$ for all $x\in X_0$, $t>0$, $r>t>0$;
\item[{\rm (ii)}] $t\mapsto T_txdt$ is Pettis integrable on $[0,+\infty[$ for all $x\in X_0$;
\item[{\rm (iii)}] $t\mapsto S_txdt$ is Pettis integrable on $[0,+\infty[$ for all $x\in X_0$.
\end{itemize}
Then $\T$ is frequently hypercyclic.
\end{theorem}

\proof We will show that $T_1$ is a frequent hypercyclic operator.  The assertion will follow
by the previous result. First observe that for any $x\in X_0$, the map $t\mapsto S_tx$ is continuous; indeed, if we fix $r>t$,
 $S_tx=T_{r-t}(S_r x)$.

 To verify that $T_1$ is frequently hypercyclic, we will follow the proof of Theorem 2.4 in
 \cite{BGE1}, by considering suitable unconditionally convergent series of integrals.

 W.l.o.g., we assume that $X_0=\{y_1, y_2, \dots\}$ is a countable set. Conditions (ii), (iii) and
 Corollary \ref{uncocon} imply that there is an increasing sequence $\{N_l\}_{l\in\N}$ in $\N$ such that,
 for all $\lambda \leq l$ and all compact sets $K\subset [N_l, +\infty[$, we have

 \begin{eqnarray}\label{2:eq1}
  & &\norm{ \int_K T_ty_\lambda dt} <\frac{1}{l2^l},\\
  \label{2:eq2}  & &\norm{ \int_K S_ty_\lambda dt} <\frac{1}{l2^l}.
 \end{eqnarray}

For every $l,\nu\in\N$, set $\rho(l, \nu)=\nu$ and apply Lemma 2.5 in \cite{BGE} to find pairwise
disjoint sets $A(l,\nu)\subseteq \N$, $l,\nu\in\N$, of positive lower density  such that,
for all $n\in A(l,\nu)$, $m\in A(k,\nu)$, with $n\not=m$,

\begin{equation}\label{2:eq3}
n\geq \nu, \quad \abs{n-m}\geq \nu + \mu.
\end{equation}

Define now

\begin{equation}\label{2:eq4}
z_n=\left\{
\begin{array}{ll}
y_l, \quad &n\in A(l,N_l)\\
0, \quad &\mbox{ otherwise}
\end{array}
\right.
\end{equation}
and set

\begin{equation}\label{2:eq5}
x:=\sum_{n\geq 1} \int_n^{n+1} S_tz_n dt.
\end{equation}

To see that this series is convergent, observe that, for each $l\in\N$,

\begin{equation}
\sum_{n\in A(l, N_l)} \int_n^{n+1}S_tz_ndt = \sum_{n\in A(l, N_l)} \int_n^{n+1} S_t y_l dt
\end{equation}
converges unconditionally by (\ref{2:eq2}). On the other hand, for every finite subset
$F\subset A(l,N_l)$, by (\ref{2:eq3})
we have that $\bigcup_{n\in F} [n.n+1]\subset [N_l, +\infty[$, hence, by (\ref{2:eq3}),

\begin{equation}
\norm{ \sum_{n\in F} \int_n^{n+1} S_ty_ldt}  \leq \frac{1}{l2^l}.
\end{equation}
Therefore we easily obtain that the series in (\ref{2:eq5}) is convergent.

Fix $l\in\N$ and $n\in A(l, N_l)$. Then

\begin{eqnarray}\label{2:eq6}
 T_{n+1}x = & & \nonumber \\
  = & \sum_{j<n} T_{n+1}  \left(\int_j^{j+1}S_tz_jdt\right) + T_{n+1} \left( \int_n^{n+1}
S_tz_ndt\right) + \sum_{j>n}T_{n+1}\left(\int_j^{j+1}S_tz_jdt\right) & \nonumber\\
  = & \sum_{j<n} \int_j^{j+1}T_{n+1-t}z_jdt + \int_n^{n+1} T_{n+1-t}y_ldt +
\sum_{j>n}\int_j^{j+1}S_{t-n-1}z_jdt  & \nonumber\\
    &= \sum_{m<n} \int_m^{m+1}T_sz_{n-m}ds + u_l +
 \sum_{m=1}^\infty\int_{m-1}^{m}S_rz_{n+m}dr, &
\end{eqnarray}

where $u_l=\int_0^1T_ty_ldt$.

We analyze the first summand in (\ref{2:eq6}):

\begin{eqnarray*}
& &\sum_{m<n}\int_m^{m+1} T_sz_{n-m}ds=\\
&=& \sum_{\lambda=1}^l \sum_{n-m\in A(\lambda, N_\lambda), m<n} \int_m^{m+1}T_sy_\lambda ds
+ \sum_{\lambda>1} \sum_{n-m\in A(\lambda, N_\lambda), m<n} \int_m^{m+1}T_sy_\lambda ds.
\end{eqnarray*}

By (\ref{2:eq3}), since $n\in A(l, N_l)$, $n-m\in A(\lambda, N_\lambda)$, we get that
$m=n-(n-m)\geq N_l+N_\lambda$.
Thus

\begin{equation}
\norm{ \sum_{m<n}\int_m^{m+1} T_sz_{n-m}ds} \leq \sum_{\lambda=1}^l \frac{1}{l2^l} +
\sum_{\lambda >l}\frac{1}{\lambda 2^\lambda} <\frac{2}{2^l}.
\end{equation}

Analogously, by (\ref{2:eq2}), we get

\begin{equation}
\norm{ \sum_{m=1}^\infty \int_{m-1}^m S_rz_{n+m}dr} <\frac{2}{2^l},
\end{equation}
which gives, for every $n\in A(l, N_l)$

\begin{equation}
\norm{ T_{n+1}x-u_l} < \frac{4}{2^l}.
\end{equation}

Since $A(l,N_l)$ has positive lower density for each $l\in\N$, we are done if we show that
$(u_l)_l$ is a dense sequence in $X$. Indeed, $u_l=Ry_l$, $l\in\N$,
where $R$ is the continuous operator defined by

\[
Rx:=\int_0^1 T_txdt.
\]

We need to prove that $R$ has dense range. First observe that $I-T_1$ has dense range.
Otherwise there would exists $\phi\in X'$, $\phi \not=0$, such that

\[
\langle \phi, x-T_1x\rangle=0 \mbox{ for all } x\in X.
\]

This implies that for every $n\in\N$ and $x\in X$:

\[
\langle \phi, x\rangle=\langle \phi,T_nx\rangle=0\ {\rm for\ all}\ x\in X.
\]

In particular, if $0<s$, then

\[
\int_n^{n+s} \langle \phi, T_ty_l\rangle>dt = \int_0^s \langle \phi, T_{u+n}y_l\rangle du=
\int_0^s \langle \phi, T_uy_l\rangle du.
\]

The left term tends to $0$, by (\ref{2:eq1}), as $n\rightarrow \infty$. Since the
right term is fixed and $s>0, l\in \N$ were arbitrary, we have

\[
\langle \phi, x\rangle =0 \mbox{ for all } x\in X,
\]
which is a contradiction.
Finally observe that if $(A, D(A))$ is the generator of $\T$, then for every $x\in D(A)$,

\[
(I-T_1)x=\int_0^1 T_t Ax = R(Ax),
\]
thus $(I-T_1)(D(A)) \subseteq R(X)$. By the density of $D(A)$ in $X$, we get that

\[
X=\overline{(I-T_1)(X)}= \overline{(I-T_1)(D(A))} \subseteq \overline{R(X)} \qquad \qed.
\]

\bigskip

\begin{corollary}\label{cor:dsw} Let $X$ be a separable complex Banach space,
$\T$ a $C_0$-semigroup with generator $A$.
 Assume that there exists a family $(f_j)_{j\in \Gamma }$ of locally bounded measurable maps
 $f_j:I_j \rightarrow X$ such that $I_j$ is an interval in $\R$,
                    $Af_j(t)=itf_j(t)$ for every $t \in I_j$, $j\in \Gamma$ and
                     $\spa \{f_j (t) \ : \ j\in \Gamma,\ t\in I_j\}$ is dense in $X$.

If either

a) $f_j \in C^2(I_j, X)$ for every $j\in \Gamma$

or

b) $X$ does not contain $c_0$  and  for every $\varphi\in X'$ and $j\in\Gamma$  the map
$\langle \varphi,f_j\rangle \in C^{1}(I_j)$,
then $\T$ is frequently hypercyclic.
\end{corollary}

\proof.

First we prove the following:
\begin{itemize}
\item[{\rm a)'}] if a) holds then there exists a family $(g_j)_{j\in \Lambda}$ of  functions
$g_j\in C^2(\R, X)$ with compact support such that
                              $Ag_j(t)=itg_j(t)$ for every $t \in \R$, $j\in \Lambda$ and
                     $\spa \{g_j (t) \ : \  j\in \Lambda,\ t\in \R\}$ is dense in $X$.
\item[{\rm b)'}] if b) holds then there exists a family $(g_j)_{j\in \Lambda}$ of  bounded
measurable functions $g_j:\R \rightarrow  X$ with compact support such that $\langle \varphi,
g_j\rangle \in C^1(\R)$ for every $\varphi \in X'$,
                              $Ag_j(t)=itg_j(t)$ for every $t \in \R$, $j\in \Lambda$ and
                     $\spa \{g_j (t) \ : \ j\in \Lambda,\ t\in \R\}$ is dense in $X$.
\end{itemize}

If $I_j=]x_j-r_j, x_j + r_j[$ is a bounded interval, consider a sequence $(\phi_n^j)_n\in C^\infty(\R)$
such that $\phi_n^j(s)=1$ if
$\abs{s-x_j}\leq r-\frac{1}{2n}$ and $\phi_n^j(s)=0$ if $\abs{s-x_j}>r_j-\frac{1}{n}$.
If we extend $f_j$ outside $I_j$ setting $f_j=0$ in $\R\setminus I_j$, we have that
$\phi_n^j f\in C^2(\R, X)$ for every $n\in\N$ if a) holds and $\langle \varphi, \phi_n^j f_j\rangle \in C^1(\R)$
for every $\varphi \in X'$ if b) holds. Moreover  $(\phi_n^j f_j)_n$ converges pointwise to $f_j$ and

\[
A(\phi_n^j(t) f_j(t))=\phi_n^j A f_j(t)= it\phi_n^j f_j(t)
\]
 for every $t \in \R$, $j\in \Gamma$.

If the interval $I_j$ is unbounded, for example $I_j=]a_j, +\infty [$, the argument runs analogously,
by considering functions $\phi_n^j\in C^\infty(\R)$ with support in $]a_j+\frac{1}{n},n[$.

It remains to show that

\begin{equation}\label{eq:span}
\overline{\spa \{\phi_n^jf_j (t) \ : \ j\in \Gamma,\ t\in I_j,\
n\in\N\}}=X.
\end{equation}

If $\varphi \in X'$ and $\langle \varphi, \phi_n^j(t) f_j(t)\rangle=0$ for every $j\in \Gamma$,
$t\in I_j$, $n\in\N$, then, by taking the limit as $n\to \infty$, we get that
$\langle \varphi, f_j(t)\rangle=0$ for every $t\in I_j$ and $j\in \Gamma$.
Then, by the assumption on the ranges of the $f_j$, we get that $\varphi=0$.

>From now on, let $\Lambda=\{ (j,n) \ : \  j\in \Gamma,\ n\in \N\}$ and,
for every $\lambda=(j,n)\in \Lambda$, set $g_\lambda=\phi_n^j f_j$.

For every $r\in\R$ and $\lambda\in \Lambda$, set

\[
\psi_{r,\lambda}:=\int_{\R} e^{-irs}g_\lambda(s)ds = \mathcal{F}(g_\lambda)(r),
\]
where $\mathcal F$ denotes the $X$-valued Fourier transform. The set
$\{\psi_{r,\lambda} \ : \  r\in\R,\ \lambda\in\Lambda\}$ is dense in
$X$. Indeed, let $\varphi \in X'$ such that for all $r\in \R$,
$\lambda\in \Lambda$,

 \[
 \langle \varphi ,\psi_{r,\lambda} \rangle = \int_{\R} e^{-irs}\langle \varphi, g_\lambda \rangle ds=0.
 \]

This means that the Fourier transform of the (scalar) function $s\mapsto \langle \varphi, g_\lambda \rangle$
vanishes on $\R$, hence, observing that $\langle \varphi, g_\lambda \rangle$ is continuous, we get that
$\langle \varphi, g_\lambda \rangle=0$   on $\R$, hence  $\varphi=0$.

For every $t>0$ set

\[
S_t\psi_{r,\lambda}:=\int_{\R} e^{-i(t+r)s} g_\lambda(s)ds = \mathcal{F}(g_\lambda)(t+r) =\psi_{t+r,\lambda}.
\]

It holds

\[
T_t\psi_{r,\lambda}=\int_{\R} e^{i(t-r)s} g_\lambda(s)ds = \mathcal{F}(g_\lambda)(-t+r)=\psi_{-t+r,\lambda},
\]
and $T_tS_t\psi_{r,\lambda}=\psi_{r,\lambda}$, $T_tS_s\psi_{r, \lambda}=S_{s-t}\psi_{r,\lambda}$
for all $\lambda\in\Lambda$, $r\in\R$, $s>t>0$.

It remains to show that $t\mapsto S_t\psi_{r,\lambda}dt$ and $t\mapsto T_t\psi_{r, \lambda}dt$
are Pettis integrable on $[0,+\infty[$ for every $r\in\R$ and $\lambda\in\Lambda$.

In the case   \textit{a)'},  $\mathcal{F}(g_\lambda)$ is Bochner integrable. Indeed,
$g_\lambda \in C^2(\R, X)$ and has compact support. Hence $g_\lambda''$ is Fourier integrable and

\[
\mathcal{F}(g_\lambda'')(r)=-r^2 \mathcal{F}( g_\lambda).
\]

Therefore $\mathcal{F}(g_\lambda)$ is  Bochner integrable on $\R$. It follows that
$t\mapsto T_t(\psi_{r,\lambda})$ and $t\mapsto S_t(\psi_{r,\lambda})$ are Bochner integrable on $[0,+\infty[$.

In the case  \textit{b)'}, we prove  that $\mathcal{F}(g_\lambda)$ is Pettis integrable on $[0,+\infty[$.
It will follow that $t\mapsto T_t(\psi_{r,\lambda})$ and $t\mapsto S_t(\psi_{r,\lambda})$ are Pettis
integrable on $[0,+\infty[$. First observe that $\mathcal{F}(g_\lambda)$ is continuous, hence  measurable.
Let $\varphi\in X'$ and consider $g(s)=\langle \varphi, g_\lambda (s)\rangle\in C^1_c(\R)$. We have:

\[
\langle \varphi, \mathcal{F}(g_\lambda)(r)\rangle =  \int_{\R}
e^{-irs}\langle \varphi, g_\lambda(s)\rangle ds= \mathcal{F}(g)(r).
\]

We have that $g'\in L^2(\R)\cap L^1(\R)$ and  $\mathcal{F}(g')(r)=ir \mathcal{F}(g)(r)\in L^2(\R)$.
Hence, for $a>0$:

\[
\int_{\abs{r}>a} \abs{\mathcal{F}(g)(r)}dr \leq \left(\int_{\abs{r}\geq a}\frac{1}{r^2}dr\right)^{\frac 1 2}
\left(\int_{\abs{r}>a} r^2\abs{\mathcal{F}(g)}^2\right)^{\frac 1 2} <+\infty.
\]

Therefore $\mathcal{F}(g)\in L^1(\R)$.  By Theorem \ref{DU}, this implies that $\mathcal{F}(g_\lambda)$ is
Pettis integrable on $[0,+\infty[$.

\qed

\begin{remark}\label{rem:dsw}
With the same argument of  \cite[Remark 2.2]{em}, one can show  that
the Desch-Schappacher-Webb criterion for chaoticity of $C_0$-semigroups
(see \cite{dsw}) implies frequent hypercyclicity.
\end{remark}

\begin{remark} There is a connection between Corollary \ref{cor:dsw} and the recent
results of S. Grivaux in \cite{gr}.
Indeed assume that one of the conditions a) or b)  (or equivalently a)' or b)') hold for a
countable family of locally bounded functions $\{f_j\}_{j\in\N}$.
Define for every $j,k\in\N$, $\theta\in [0,2\pi[$

\[
E_{j,k}(e^{i\theta})= f_j(\theta+2k\pi).
\]

The family $\{E_{j,k} \ : \ j,k\in\N\}$ is a countable family of
bounded $\mathbb{T}$-eigenvectors fields for the operator $T_1$,
where $\mathbb{T}=\{\lambda\in\C \ : \  |\lambda|=1\}$, such that
$\spa \{E_{j,k} (\lambda) \ : \ \lambda\in \mathbb{T},\ j,k\in\N\}$
is dense in $X$. Actually it holds that  $\spa \{E_{j,k} (\lambda) \
: \  \lambda\in \mathbb{T}\setminus D,\ j,k\in\N\}$ is dense in $X$
for every countable subset of $\T$. Indeed, if $D=\{ e^{i\theta_n} \
: \  n\in\N\}$, with $\theta_n\in [0,2\pi[$, then

\[
\spa \{E_{j,k} (\lambda) \ : \  \lambda\in \mathbb{T}\setminus D\,
j,k\in\N\}= \spa \{f_{j}(s) \ : \ s\in \R \setminus \{
\theta_n+2k\pi \ : \ k\in\N\}, \ j\in\N\}.
\]

By the weak continuity of the vector field $f_j$, we get that the right hand set is dense in $X$.

Thus, by \cite[Proposition 4.2]{gr} $T_1$ has  perfectly spanning
eigenvectors associated to unimodular eigenvectors, i.e. there
exists a probability measure $\sigma$ on the unit circle
$\mathbb{T}$ such that for every $\sigma$-measurable subset of $A$
with $\sigma(A)=1$, $\spa \{\ker(T-\lambda) \ : \ \lambda\in\ A\}$
is dense in $X$.
\end{remark}

\begin{example}
Consider  the linear perturbation of the one-dimensional
Ornstein-Uhlenbeck operator

\[
{\mathcal{A}_\alpha}u=u''+bxu'+\alpha u,
\]
     where $\alpha\in \R$, with domain

\[
  D(\mathcal{A}_\alpha)=\left\{\, u\in L^2(\R)\cap
  W^{2,2}_{{\rm loc}}(\R) \ : \  \mathcal{A}_\alpha u\in L^2(\R)\,\right\}.
\]

In \cite{cm}, it was proved that if $\alpha>b/2>0$, then the semigroup generated by
$\mathcal{A}_\alpha$ in $L^2(\R)$ is chaotic. Actually the semigroup is frequently hypercyclic.

Indeed,  for every $\mu\in \C$,  with $\Re \mu < -\frac b 2 +\alpha$
the functions $u^1_\mu$ and $u^2_{\mu}$,  whose Fourier transforms are

\[
    \widehat{u^1_\mu} (\xi)=e^{-\xi^2/2b}\xi |\xi|^{-(2+(\mu-\alpha)/b)}, \qquad
    \widehat{u^2_\mu} (\xi)=e^{-\xi^2/2b}
    |\xi|^{-(1+(\mu-\alpha)/b)},
\]
    are eigenfunctions of $\mathcal{A}_\alpha$ (see \cite{cm,me}).

 For every $s\in \R$, consider the functions
        $f_1(s)=u^1_{is}$ and $f_2(s)=u^2_{is}$. For every $\phi \in X' =L^2(\R)$ and $j=1,2$, by Parseval equality, we have

\[
\langle \phi, f_j (s)\rangle= \int_\R \phi(x)u_{is}^j(x) dx= \int_\R \widehat\phi(x)\widehat u_{is}^j(x) dx \qquad s\in \R.
\]

It is immediate to verify that  $\langle\phi, f_j\rangle\in C^1(\R)$, by Lebesgue's theorem.

    The argument of \cite{cm} shows that $\spa \{f_i(s) \ : \  i=1,2, s\in \R\}$ is dense in $L^2(\R)$.
Therefore  the semigroup is frequently hypercyclic by  Corollary
\ref{cor:dsw}.
\end{example}

We will see that the Frequent Hypercyclicity Criterion for semigroups implies chaos for each single operator
of the semigroup. It is interesting to observe that this is in general stronger than the
chaoticity of the semigroup since, by recent results of Bayart and Berm\'{u}dez \cite{bayart_bermudez2009semigroups},
there are chaotic $C_0$-semigroups $\T$ such that no single operator $T_t$ is chaotic, and
chaotic $C_0$-semigroups $\T$ containing non-chaotic operators $T_{t_0}$, $t_0>0$, and at the same
time chaotic $T_{t_1}$ for some $t_1>0$.

\begin{proposition}\label{pr:fhcrit_chaos}
Let $X$ be a separable Banach space and let $\T$ be a $C_0$ semigroup on $X$ that satisfies the Frequent Hypercyclicity
Criterion of Theorem~\ref{fhcriterion}. Then the operator $T_{t_0}$ is chaotic for every $t_0>0$.
\end{proposition}

\proof Given $t_0>0$, we know that $T_{t_0}$ is frequently hypercyclic, thus hypercyclic \cite{conejero_muller_peris2007hypercyclic}.
Given $x\in X$ and $\varepsilon >0$ we want to find a $T_{t_0}$-periodic point $z\in X$ such that $\norm{x-z}<\varepsilon$.
Indeed, let $y\in X_0$ such that $\norm{x-y}<\varepsilon $. By continuity, we fix $\delta >0$ such that

\[
\norm{x-\delta^{-1}\int_0^\delta T_syds}<\varepsilon .
\]

Now, let $n\in \N$ big enough so that, by Corollary~\ref{uncocon} and for $t:=nt_0$, we have that

\[
z:= \delta^{-1} \left[ \sum_{k\geq 1}  \int_0^\delta S_{kt-s}yds   +    \int_0^\delta T_syds   +
\sum_{k\geq 1}  \int_0^\delta T_{kt+s}yds \right]
\]
satisfies $\norm{x-z}<\varepsilon$. Finally, observe that the hypothesis of the Frequent Hypercyclicity Criterion
and continuity give $T^n_{t_0}z=T_tz=z$.  \qed

Finally we point out the connection between the Frequent Hypercyclicity Criterion for semigroups and the
Frequent Hypercyclicity Criterion for sequences of operators.

\begin{proposition}\label{pr:fhcritop}
Let $X$ be a separable Banach space and let $\T$ be a $C_0$ semigroup on $X$. If
there  exists $X_0\subseteq X$, dense subset, with $T_t(X_0)\subseteq X_0$ for every
$t>0$, and maps $S_t:X_0\rightarrow X_0$, $t>0$, satisfying:
\begin{itemize}
\item[{\rm (i)}] $T_tS_tx=x$, $S_rT_tx=T_tS_rx=S_{r-t}x$ for all $x\in X_0$, $t>0$, $r>t>0$;
\item[{\rm (ii)}] $t\mapsto T_txdt$ is Pettis integrable on $[0,+\infty[$ for all $x\in X_0$;
\item[{\rm (iii)}] $t\mapsto S_txdt$ is Pettis integrable on $[0,+\infty[$ for all $x\in X_0$.
\end{itemize}
then the operator $T_{t_0}$ satisfies the Frequent Hypercyclicity criterion for operators for every $t_0>0$.
\end{proposition}

\proof By sake of simplicity, let $t_0=1$.  First observe that
$S_nx=S_nT_1S_1x=S_{n-1}S_1x=S_{n-2}S_1^2=...=S_1^nx$ for every $x\in X_0$.  For every $x\in X_0$, set
$y=\int_0^1 T_tx dt$. Then the series

\[
\sum_{n=1}^\infty T_n y = \sum_{n=1}^\infty \int_n^{n+1}T_tx dt
\]
is unconditionally convergent by Proposition \ref{ucPettis}.
Analogously, since

\[
\int_{n-1}^{n}S_sxds=\int_0^1 S_{n-1+s}x ds =\int_0^1 S_{n-u}xds =S_n \int_0^1T_sxds=S_n y,
\]
we get that the series $\sum_{n=1}^\infty S_ny$ is unconditionally
convergent. Finally we observe that, by  the same argument used in
the proof of Theorem \ref{fhcriterion}, the set $\{ \int_0^1 T_tx \
: \ x\in X_0\}$ is dense. \qed

Given a $C_0$-semigroup $\T$, it seems a reasonable guess that if
$T_{t_0}$ satisfies the Frequent Hypercyclicity Criterion for
operators of Proposition~\ref{pr:frhyop} for some $t_0>0$, then the
semigroup $\T$ should satisfy the Frequent Hypercyclicity Criterion
for semigroups of Theorem~\ref{fhcriterion}. Unfortunately, we do
not  know whether  this holds in general.

\section{Translation semigroups}

An \textit{admissible weight} function on $[0,+\infty[$ is a
measurable function $\rho : [0,+\infty[ \rightarrow \R$ satisfying
the following conditions:

(i) $\rho(t) > 0$ for all $t\in [0,+\infty[ $;

(ii) there exist constants $M \geq  1$ and $\omega \in \R$ such that
$\rho(\tau)\leq Me^{\omega t} \rho(\tau+t)$ for all
$\tau\in [0,+\infty[ $ and all $t > 0$.

We recall the following useful results:

\begin{lemma}\label{le:dis}
If $\rho$ is an admissible weight, then for every $l>0$ there exist $A, B>0$
such that for every $\sigma\in [0,+\infty[ $ and for every $t\in [\sigma, \sigma + l]$, it holds
\[
A\rho(\sigma) \leq \rho(t)\leq B\rho(\sigma +l).
\]
\end{lemma}

\begin{lemma}\label{le:int}
Let $\rho:[0,+\infty[ \rightarrow \R^+$ be an admissible weight.
\begin{enumerate}
\item The following conditions are equivalent:

\begin{itemize}
\item[{\rm (i)}] For all $b\geq 0$ the series $\sum_{k=1}^\infty \rho (b+k)$  is convergent,
\item[{\rm (ii)}] For all $b\geq 0$ there exists $P>0$ such that the series
$\sum_{k=1}^\infty \rho (b+kP)$ is convergent,
\item[{\rm (iii)}] There exists $P>0$ such that the series
$\sum_{k=1}^\infty \rho (kP)$ is convergent,
\item[{\rm (iv)}] The series $\sum_{k=1}^\infty \rho (k)$ is convergent,
\item[{\rm (v)}] $\int_0^{+\infty} \rho(s)ds<+\infty$,
\item[{\rm (vi)}] there exists   $D\subseteq \N$ with bounded gaps (i.e. there exists $M>0$
such that $D\cap [n,n+M]\not=\emptyset$ for every $n\in\N$) such that $\sum_{k\in D} \rho(k)$ is convergent.
\end{itemize}

\item $\rho$ is bounded if and only if there exists $D\subseteq \N$
with bounded gaps such that $\rho$ is bounded on $D$.
\end{enumerate}
\end{lemma}

\proof  1. The implications (i)$\Rightarrow$ (ii) and (ii) $\Rightarrow$ (iii) are obvious.
If (iii) holds,  we fix $n\in \N$ with $n>P$, and $A,B>0$ satisfying
the inequalities of Lemma~\ref{le:dis} for $l=P$. For all $k \in \N$ there exists
$m_k\in \N$ such that $kn\in [(m_k-1)P,m_kP]$.
Since $n>P$ then $m_{k+1}>m_k$ for all $k\in \N$. Thus,

\[
\sum_{k=1}^\infty \rho (kn)\leq B\sum_{k=1}^\infty \rho (m_kP)\leq B\sum_{k=1}^\infty \rho (kP)<+\infty .
\]

Analogously, $\sum_{k=1}^\infty \rho (kn+j)\leq B \sum_{k=1}^\infty
\rho (kP)$, $j=1,\dots ,n-1$, and we obtain that $\sum_{k=1}^\infty
\rho (k)<+\infty$.

If (iv) holds, then

\[
\int_0^{+\infty} \rho (s)ds =\sum_{k=1}^\infty \int_{k-1}^k \rho(s)ds\leq B\sum_{k=1}^\infty \rho (k)< +\infty
\]
by an application of Lemma \ref{le:dis} for $l=1$, which yields (v).

If $\int_0^{+\infty} \rho(s)ds<+\infty$ and $b\geq 0$ is arbitrary, then

\[
\sum_{k=1}^\infty \rho (b+k) \leq A^{-1} \sum_{k=1}^\infty \int_{b+k}^{b+k+1}
\rho (s)ds =A^{-1} \int_{b+1}^{+\infty} \rho (s)ds<+\infty ,
\]
where we applied again Lemma \ref{le:dis} for $l=1$. This concludes that (v) implies (i).

Finally, clearly (iii) implies (vi). We prove that (vi) implies (iii). Assume w.l.o.g. that
$M\in\N$. Then for every $k\in\N$ there exists $h_k\in D\cap [Mk,M(k+1)]$. If $A,B>0$ satisfy
the inequalities of Lemma~\ref{le:dis} for $l=M$, we have that

\[
A\rho(Mk)\leq \rho(h_k) \leq B\rho(M(k+1)), \qquad k\in\N.
\]

Thus

\[
\sum_{k=1}^\infty \rho(Mk) \leq A^{-1} \sum_{k=1}^\infty \rho(h_k) \leq A^{-1} \sum_{h\in D} \rho(h)<+\infty.
\]

2.  Assume that $\rho(h)\leq K$ for every $h\in D$, where  $D\subseteq \N$ with bounded gaps. Hence,
there is $M\in\N$ such that $[Mn,Mn+M]\cap D\not=\emptyset$ for every $n\in\N$.
Choose $h_n\in [Mn, Mn+M]\cap D$. By Lemma \ref{le:dis}, there exists $A_M, B_M>0$ such that

\[
A_M\rho(Mn)\leq \rho(h_n)\leq B_M\rho(Mn+M), \qquad n\in\N,
\]
hence $\rho(Mn)\leq A_M^{-1}K$ for every $n\in\N$. On the other hand, for every $s\geq0$ there
exists $\overline  n$ such that $x\in [M\overline n, M\overline n+M]$,and therefore

\[
\rho(s) \leq B_M\rho(M\overline n +M) \leq KA_M^{-1} B_M. \qquad \qed
\]

\medskip

 Next we consider the following function spaces:

\[
L_p^\rho([0,+\infty [)=\{ u: [0,+\infty [ \rightarrow \R \ : \ u
\mbox{ is measurable and } \norm{ u}_p<\infty\},
\]
where $\norm{ u}_p=\left(\int_0^\infty \abs{u(t)}^p\rho(t) dt\right)^{\frac 1 p}$, and

\[
C_0^\rho([0,+\infty [)=\{u:[0,+\infty [ \rightarrow \R \ : \ u
\mbox{ is continuous and } \lim_{x\to\infty}u(x)\rho(x)=0\},
\]
with $\norm{u}_\infty = \sup_{t\in [0,+\infty [}\abs{u(t)}\rho (t)$.

If $X$ is any of the spaces above, the translation semigroup $\T$ is defined as usual by
\[ T_tf(x)=f(x+t),\qquad t\geq 0,\ x\in I,\]
and is a $C_0$-semigroup (see e.g. \cite{dsw}).

Hypercyclic and chaotic  translation semigroups have been characterized in \cite{dsw,delaubenfels_emamirad2001chaos,ma1}.
If   $X$ is one of the spaces $L_p^\rho([0,+\infty [)$ or $C_0^\rho([0,+\infty [)$ with
an admissible weight function $\rho$,  the translation semigroup $\T$  on $X$ is hypercyclic if
and only if
$\liminf_{t\to +\infty} \rho(t)=0$.

If   $X=C_0^\rho([0,+\infty [)$, then the
 translation semigroup $\T$ on $X$ is chaotic if and only if $\lim_{x\to+\infty}\rho(x)=0$.

If $X=L_p^\rho([0,+\infty [)$, $\T$ is chaotic if and only if any of the conditions of lemma \ref{le:int} is satisfied.

\begin{proposition} Let $\rho$ be an admissible weight on $[0,+\infty[$, $X=L^p_\rho([0,+\infty[)$, $1\leq p<+\infty$ and
$\T$ the translation semigroup on $X$. Then $\T$ is  chaotic if and only if
it satisfies the Frequent Hypercyclicity Criterion for semigroups.
\end{proposition}

\proof  Let $X=L^p_\rho([0,+\infty[)$. If $\T$ is chaotic, then $\int_0^{+\infty}\rho(s)ds$ is finite.
Let $X_0$ be the space
generated by the characteristic functions of bounded intervals of $[0,+\infty[$. $X_0$ is dense
in $L^p_\rho([0,+\infty[)$.  For every $t>0$ and $f\in X_0$ set

\[
S_t f (s)=\left\{\begin{array}{ll} f(s-t) \qquad &s\geq t\\
0 \qquad &s\in [0,t[\end{array}\right.
\]

Observe that  $T_tS_tf=f$ and  $T_tS_rf=S_{r-t}f$ for all $f\in
X_0$, $t>0$, $r>t>0$. Moreover $\int_{\R^+}\norm{T_tf}dt$ converges
for all $f\in X_0$, because of the compact support of $f$, hence
$\int_{\R^+}T_tf$ is Pettis integrable. On the other hand, consider
$f=\chi_{[a,b]}$, with $0\leq a<b$. If $p=1$, we have

\[
||S_tf||= \int_t^{t+b}\rho(s)ds=\int_0^b  \rho(s+t)ds \leq  bB\rho (t+b)
\]
where $B$ is a  positive costants such that for every $s\in [0,b]$ and for every $t\geq 0$

\[
\rho(s+t) \leq B\rho(t+b).
\]

  Since  $\int_0^\infty \rho(t+b) dt$ is finite, we get that $t\mapsto S_tf$ is Pettis integrable.
Let $p>1$ and let $\phi\in L^{p'}_\rho([0,+\infty[)$, where $\frac 1 p + \frac{1}{p'}=1$.
To prove that $t\mapsto S_tf$ is Pettis integrable, by Theorem \ref{DU}, we have to show that
$t\mapsto \langle \phi, S_tf\rangle \in L^1([0,+\infty[)$.
It holds that

\[
\langle \phi, S_t f\rangle = \int_t^{+\infty} f(s-t)\rho(s)ds=\int_0^{+\infty} f(u) \rho(t+u) du.
\]

A straightforward application of Tonelli and Fubini theorems (as for
the proof of the integrability of the convolution) gives the
assertion. \qed

\begin{proposition} Let  $\rho$ be an admissible weight on $[0,+\infty[$ and
$(T_t)_{t\geq 0}$ the translation
semigroup on $C_0^\rho([0,+\infty[)$. If $\int_0^{+\infty}
\rho(s)ds<+\infty$, then $\T$  satisfies the frequent hypercyclicity
criterion for semigroups.
\end{proposition}

\proof Let $X_0$ be the space of continuous functions on
$[0,+\infty[$ with compact support. For every $t>0$ and each $f\in
X_0$ set

\[
S_t f (s)=\left\{\begin{array}{ll} f(s-t) \qquad &s\geq t\\
f(0) \qquad &s\in [0,t[\end{array}\right.
\]

Observe that  $T_tS_tf=f$ and  $T_tS_rf=S_{r-t}f$ for all $f\in
X_0$, $t>0$, $r>t>0$. Moreover, because of the compact support of
$f$, $T_t f=0$ for $t$ big enough. Hence
$\int_{0}^{+\infty}\norm{T_tf}dt$ converges for all $f\in X_0$. On
the other hand, if $f\in X_0$ and $\supp f\subseteq [0,b]$, then

\[
\norm{S_tf}= \sup_{s\in [t,t+b]}|f(s-t)|\rho(s) =\sup_{s\in [0,b]}
|f(s)| \rho(s+t) \leq  BA^{-1}\rho(0)^{-1}\rho (t+b) \norm{f},
\]
where $A,B$ are  positive constants such that for every $s\in [0,b]$
and for every $t\geq 0$

\[
A\rho(t)\leq \rho(s+t) \leq B\rho(t+b).
\]

  Since  $\int_0^{+\infty} \rho(t+b)dt$ is finite, we get that
  $\int_0^{+\infty} S_tfdt$ converges unconditionally.
\qed

\begin{remark}
It should be observed that, for the translation semigroup $\T$ on
$L_p^\rho([0,+\infty[)$, it holds that $\T$ is chaotic if and only
if  every operator $T_t$  satisfies the Frequent Hypercyclicity
Criterion for operators by Proposition \ref{pr:fhcritop}.
\end{remark}

In \cite{grosse-erdmann_peris2005frequently} the authors obtain a necessary condition for frequent hypercyclicity
of unilateral weighted shifts on $\ell^p$. Inspired by their condition, we obtain an analogous one for
translation semigroups.

\begin{proposition} Let $\rho$ be an admissible weight on $[0,+\infty[$, $\T$ the translation semigroup in
$L_p^\rho([0,+\infty[)$. If $\T$ is frequently hypercyclic, then for every $\varepsilon>0$ there exists
a sequence $(n_k)_k$ in $\N$ with positive lower density such that

\begin{equation}\label{eq:neccon}
\sum_{k> i} \rho(n_k-n_i)<\varepsilon.
\end{equation}

Moreover,  $\rho$ is bounded.
\end{proposition}

\proof By the results in \cite{conejero_muller_peris2007hypercyclic}, the operator
$T_1$ is frequently hypercyclic. Let  $f\in L_p^\rho([0,+\infty[)$ be frequently hypercyclic
for $T_1$. If $u$ is the characteristic function of $[0,1]$, $\rho_0=(\int_0^1\rho(s)ds)^{-\frac 1 p}$,  and
$0<\eta<1$, there exists a sequence $(n_k)_k$ in $\N$ with positive lower density such that
$\norm{T_{n_k}f-\rho_0u}<\eta$, $k\in \N$.
This implies that

\[
(\int_0^1 |f(s+n_k)-\rho_0|^p\rho(s)ds )^{\frac 1 p}< \eta,
\]
hence

\begin{eqnarray*}
 & &1=\rho_0 (\int_0^1 \rho(s)ds)^{\frac 1 p} \leq \\
&\leq& (\int_0^1 \abs{f(s+n_k)-\rho_0}^p\rho(s)ds )^{\frac 1 p} +
(\int_0^1 \abs{f(s+n_k)}^p\rho(s)ds )^{\frac 1 p}<\\
&<& \eta + (\int_0^1 \abs{f(s+n_k)}^p\rho(s)ds )^{\frac 1 p}.
\end{eqnarray*}

So we get that

\[
\int_0^1 \abs{f(s+n_k)}^p\rho(s)ds > (1-\eta)^p.
\]

By lemma \ref{le:dis}, there exist $A,B>0$ such that $A\rho(\tau)\leq \rho(s)\leq B\rho(\tau+1)$
for every $\tau\geq 0$ and $s\in [\tau,\tau+1]$. Hence

\[
\int_0^1 |f(s+n_k)|^pds\geq \frac{1}{B\rho(1)}\int_0^1 |f(s+n_k)|^p\rho(s)ds
>\frac{1}{B\rho(1)}(1-\eta)^p, \ \ k\in \N.
\]

For all $i\in\N$ we have

\begin{eqnarray*}
 & &\eta^p>\int_0^{+\infty}|f(s+n_i)-\rho_0 u(s)|^p\rho(s)ds
 >\int_1^{+\infty}|f(s+n_i)|^p\rho(s)ds \geq\\
&\geq &\int_{n_i+1}^\infty |f(s)|^p\rho(s-n_i)ds=\sum_{j=1}^\infty
\int_{n_i+j}^{n_i+j+1}|f(s)|^p\rho(s-n_i) ds \geq \\
&\geq& \sum_{j> i}^\infty \int_{n_j}^{n_j+1}|f(s)|^p\rho(s-n_i) ds=\\
&= &\sum_{j> i}\int_0^1 |f(s+n_j)|^p\rho(s+n_j-n_i)ds\geq \\
&\geq&A\sum_{j> i}\rho(n_j-n_i) \int_0^1 |f(s+n_j)|^pds >
\frac{A(1-\eta)^p}{B\rho(1)} \sum_{j> i} \rho(n_j-n_i).
\end{eqnarray*}

Hence

\[
\sum_{j> i} \rho(n_j-n_i) < \frac{B\rho(1)}{A(1-\eta)^p}\eta^p .
\]

Since $\eta\in ]0,1[$ was arbitrary, we get the desired inequality.
As a consequence of (\ref{eq:neccon}), since for every $A\subset\N$ with positive lower
density the difference set $A-A$ has bounded  gaps (see \cite{ST}), we have that there
exists $D\subseteq \N$ with bounded gaps such that $\rho(h)<1$ for every $h\in D$.
Hence $\rho$ is bounded by Lemma \ref{le:int}(2). \qed

\begin{example} Let $\phi:\R_+ \rightarrow \R$ be a $C^1$-function with  derivative
bounded by $\omega>0$ and such that

\[
\limsup_{t\to+\infty} \phi(s)=+\infty,\ \liminf_{t\to-\infty} \phi(s)=-\infty.
\]
(For example, consider a $C_1$ function such that $\phi(s)=s\sin(\log s)$ if $s\geq 1$.)
Set $\rho=e^{-\phi}$. Clearly $\rho>0$ and if $t,\tau>0$ we have

\[
\frac{\rho(\tau)}{\rho(t+\tau)}= e^{-\int_\tau^{t+\tau}\phi'(s)ds} \leq e^{\omega \tau}.
\]

Hence $\rho$ is an admissible weight. The translation semigroup in $L_p^\rho([0,+\infty[)$
is hypercyclic, since $\liminf_{s\to\infty}\rho(s)=0$, but it is not frequently hypercyclic,
since $\rho$ is unbounded.
\end{example}

With a similar proof we get a necessary condition  for frequently hypercyclic
translation semigroups in weighted spaces of continuous functions.

\begin{proposition} Let $\rho$ be an admissible weight on $[0,+\infty[$, and $\T$ the translation
semigroup in $C_0^\rho([0,+\infty[)$. If $\T$ is frequently hypercyclic, then for every
$\varepsilon>0$ there exists a sequence $(n_k)_k$ in $\N$ with positive lower density such that
$\rho(n_k-n_i)<\varepsilon$ for every $k> i$.
\end{proposition}

\section{Appendix}

We recall in this appendix the main definitions and results about Pettis integrability.
Let $X$ be a Banach space and $(\Omega,\mu)$  a $\sigma$-finite measure space. A function
$f:\Omega\rightarrow X$ is said to be weakly $\mu$-measurable if the scalar function $\varphi\circ f$
is $\mu$-measurable for every $\varphi\in X'$, where  $X'$  denotes the topological dual of $X$;
$f$ is said to be $\mu$-measurable if there exists a sequence $(f_n)_n$ of simple functions
 such that $\lim_{n\to\infty}\abs{f_n-f}=0$ $\mu$-a.e.

\begin{lemma}[Dunford] If $f$ is weakly $\mu$-measurable and $\varphi\circ f\in L^1(\Omega,\mu)$
for every $\varphi\in X'$, then for every measurable $E\subseteq \Omega$ there exists $x_E\in X''$ such that

\[
x_E(\varphi)=\int_E \varphi\circ f d\mu,
\]
for every $\varphi\in X'$.
\end{lemma}

\begin{definition} If $f:\Omega \rightarrow X$ is weakly $\mu$-measurable and
and $\varphi\circ f\in L^1(\Omega,\mu)$ for every $\varphi\in X'$, then $f$ is called Dunford
integrable. The Dunford integral of $f$ over a measurable $E\subseteq \Omega$ is defined by
the element $x_E\in X''$ such that

\[
x_E(\varphi)=\int_E \varphi\circ f d\mu,
\]
for every $\varphi\in X'$.

In the case that $x_E\in X$ for every measurable $E$, then $f$ is called Pettis integrable and
$x_E$ is called the Pettis integral of $f$ over $E$ and will be denoted by
$(P)-\int_E f d\mu$.
\end{definition}

Clearly the Dunford and Pettis integrals coincide if $X$ is a reflexive space.

Moreover, if $\norm{f}$ is integrable on $\Omega$ (i.e. $f$ is Bochner integrable on $\Omega$),
then $f$ is Pettis integrable on $X$.

\begin{theorem}[Pettis]\label{ucPettis} If $f$ is Pettis integrable, then for every sequence
$(E_n)_n$ of disjoint measurable sets in $\Omega$

\[
\int_{\bigcup_{n\in\N} E_n}fd\mu =\sum_{n\in\N} \int_{E_n} fd\mu,
\]
where the series converges unconditionally.
\end{theorem}

\begin{corollary}\label{uncocon}
If $f:[0,+\infty[\rightarrow X$ is   Pettis integrable on $[0,+\infty[$, then
for every $\varepsilon>0$ there exists $N>0$ such that for every compact set $K\subset [N,+\infty[$

\[
\norm{\int_K f(t)dt}<\varepsilon.
\]
\end{corollary}

\proof  Assume that there exists $\varepsilon >0$ such that for every $n\in\N$ there exists
a compact set $K_n\subseteq [n, +\infty[$ such that

\[
\norm{\int_{K_n}f(s)ds} >\varepsilon.
\]

It is easy to find a sequence $(k_n)_n$ of natural numbers such that the sets
$K_{k_n}$ are mutually disjoint. Then

 \[
 \int_{\bigcup_n K_{k_n}} f(s)ds =\sum_{n=1}^\infty \int_{K_n}f(s)ds,
 \]
 hence $\lim_{n\to\infty} \int_{K_n}f(s)ds=0$, a contradiction. \qed

\begin{theorem}\label{DU} If the Banach space $X$ does not contain $c_0$ and  $(\Omega,\mu)$
is $\sigma$-finite measure space, then a measurable Dunford integrable function
$f:\Omega\rightarrow X$ is Pettis integrable.
\end{theorem}

The proofs of all these results  can be found in \cite{DU} for the case of
finite measure space, but they easily extend to $\sigma$-finite measure spaces. In particular,
the proof of the deep Theorem \ref{DU} follows analogously to the finite measure space case (\cite{DU},  Theorem 7, p.54) taking into account the following decomposition theorem
due to J. K. Brooks (see \cite{Br}, 2. Theorem 1).

\begin{theorem} Let $(\Omega, \mu)$ be a $\sigma$-finite measure space. If $f:\Omega\rightarrow X$
is a measurable weakly integrable function. Then $f$ can be represented in the form $f=g+h$ $\mu$-a.e.
where $g$ is a bounded Bochner integrable function and $h$ assumes at most a countable number of values in $X$.
\end{theorem}

\section{Acknowledgements}
The authors thank A.~Albanese for a suggestion about the end of the proof of Theorem \ref{fhcriterion}. The
research of the second author was partially supported by the MICINN  and FEDER Projects
MTM2007-64222 and MTM2010-14909,  and by Generalitat Valenciana Project
PROMETEO/2008/101.

\end{document}